\documentclass[10pt,reqno]{amsart}
\usepackage[utf8]{inputenc}
\usepackage[english]{babel}
\usepackage[a4paper,margin=3cm]{geometry}
\usepackage{amssymb,color,mathtools,hyperref,fourier}
\numberwithin{equation}{section}
\newtheorem{theorem}{Theorem}[section]
\newtheorem{lemma}[theorem]{Lemma}
\newtheorem{remark}[theorem]{Remark}
\title[Convergence of the spectral radius of a random matrix through its characteristic polynomial]{Convergence of the spectral radius of a random matrix\\through its characteristic polynomial}%
\date{Autumn 2020, revised Summer 2021. Preprint, compiled \today.}%
\author[Bordenave]{Charles Bordenave}%
\address{Institut de Mathématiques de Marseille; CNRS; %
  Aix-Marseille Université, Marseille, France.}
\email{charles.bordenave@univ-amu.fr}%
\urladdr{http://www.i2m.univ-amu.fr/perso/charles.bordenave/}%
\author[Chafaï]{Djalil Chafaï}%
\address{CEREMADE, CNRS, UMR 7534, %
  Université Paris-Dauphine, PSL University, 75016 Paris, France.}%
\email{djalil@chafai.net} \urladdr{http://djalil.chafai.net/}%
\author[García-Zelada]{David García-Zelada}%
\address{Institut de Mathématiques de Marseille; CNRS; %
  Aix-Marseille Université, Marseille, France.}%
\email{david.garcia-zelada@univ-amu.fr}%
\urladdr{https://davidgarciaz.wixsite.com/math}%
\thanks{CB and DGZ are supported by the grant ANR-16-CE40-0024}
\keywords{Random Matrix; Spectral Radius; Gaussian Analytic Function; Central Limit Theorem; Combinatorics; Digraph; Circular Law}
\subjclass[2010]{%
Primary: 30C15, % Zeros of polynomials, rational functions, and other analytic functions of one complex variable (e.g.,
60B20; % Random matrices (probabilistic aspects) 
Secondary: 60F05 % Central limit and other weak theorems
}
\begin{document}
\begin{abstract}
  Consider a square random matrix with independent and identically distributed
  entries of mean zero and unit variance. We show that as the dimension tends
  to infinity, the spectral radius is equivalent to the square root of the
  dimension in probability. This result can also be seen as the convergence of
  the support in the circular law theorem under optimal moment conditions. In
  the proof we establish the convergence in law of the reciprocal
  characteristic polynomial to a random analytic function outside the unit
  disc, related to a hyperbolic Gaussian analytic function. The proof is short
  and differs from the usual approaches for the spectral radius. It relies on
  a tightness argument and a joint central limit phenomenon for traces of
  fixed powers.
\end{abstract}
\maketitle
%
%{\footnotesize\tableofcontents}
%
\section{Introduction and main results}

Let $\{a_{jk}\}_{j,k\geq1}$ be independent and identically distributed complex
random variables with mean zero and unit variance, namely
$\mathbb{E}[a_{11}]=0$ and $\mathbb E[|a_{11}|^2] = 1$. For all $n\geq1$, let
\begin{equation}\label{eq:An}
  A_n={(a_{jk})}_{1 \leq j,k \leq n}.
\end{equation}
We call it a Girko matrix \cite{MR3808330}. When $a_{11}$ is Gaussian with
independent and identically distributed real and imaginary parts then $A_n$
has density proportional to $\mathrm{e}^{-\mathrm{Tr}(AA^*)}$ and belongs to
the complex Ginibre ensemble \cite{MR173726}. We are interested in the matrix
$\frac{1}{\sqrt{n}}A_n$ for which each row and each column has a unit mean
squared Euclidean norm. Its characteristic polynomial at point
$z\in\mathbb{C}$ is
\begin{equation}\label{eq:pn}
  p_n(z)=\det\Bigr(z-\frac{A_n}{\sqrt{n}}\Bigr)
\end{equation}
where $z$ stands for $z$ times the identity matrix. The $n$ roots of $p_n$ in
$\mathbb{C}$ are the eigenvalues of $\frac{1}{\sqrt{n}}A_n$. They form a
multiset $\Lambda_n$ which is the spectrum of $\frac{1}{\sqrt{n}}A_n$. The
spectral radius of $\frac{1}{\sqrt{n}}A_n$ is defined by
\begin{equation}\label{eq:rho}
  \rho_n=\max_{\lambda\in\Lambda_n}|\lambda|.
\end{equation}
The circular law theorem states that the empirical measure of the elements of
$\Lambda_n$ tends weakly as $n\to\infty$ to the uniform distribution on the
closed unit disc: almost surely, for every nice Borel set
$B\subset\mathbb{C}$,
\begin{equation}\label{eq:cirlaw}
  \lim_{n\to\infty}\frac{\mathrm{card}(B\cap\Lambda_n)}{n}
  =\frac{\mathrm{area}(B\cap\overline{\mathbb{D}})}{\pi},
\end{equation}
where ``$\mathrm{area}$'' stands for the Lebesgue measure on $\mathbb C$, and
where $\overline{\mathbb{D}}=\{z\in\mathbb{C}:|z|\leq1\}$ is the closed unit
disc, see \cite{MR773436,MR3808330,MR2722794,MR2908617}. The circular law
\eqref{eq:cirlaw}, which involves weak convergence, does not provide the
convergence of the spectral radius, it gives only that almost surely
\begin{equation}\label{eq:liminf}
  \varliminf_{n\to\infty}\rho_n\geq1.
\end{equation}
Theorem \ref{th:p} provides the convergence of the spectral radius, without
extra assumptions on the entries. This result was conjectured in
\cite{MR3813992}, and improves over
\cite{MR866352,MR841088,MR863545,MR3813992,basak2019outliers}. The moments
assumptions are optimal, and the $\frac{1}{\sqrt{n}}$ scaling is no longer
adequate for entries of infinite variance, see for instance \cite{MR2837123}.

We have $\rho_n\leq\sigma_n$ where $\sigma_n$ is the operator norm of
$\frac{1}{\sqrt{n}}A_n$, its largest singular value. It is known that the
condition $\mathbb{E}[|a_{11}|^4]<\infty$ is necessary and sufficient for the
convergence of $\sigma_n$ as $n\to\infty$, see \cite{MR2567175}. A stricking
aspect of the spectral radius is that it converges without any extra moment
condition.

\begin{theorem}[Spectral radius]\label{th:p}
  We have $\lim_{n\to\infty}\rho_n=1$ in probability, in the sense that for
  all $\varepsilon>0$,
  \[
    \lim_{n\to\infty}\mathbb{P}(|\rho_n-1|\geq\varepsilon)=0.
  \]
\end{theorem}

The proof of Theorem \ref{th:p} is given in Section \ref{se:pr:th:p}. It
relies on Theorem \ref{th:q} below, which is of independent interest. It does
not involve any Hermitization or norms of powers in the spirit of Gelfand's
spectral radius formula. The idea is to show that on
$\mathbb{C}\cup\{\infty\}\setminus\overline{\mathbb{D}}$, the polynomial
$z^{-n}p_n(z)$ tends as $n\to\infty$ to a random analytic function which does
not vanish. The first step for mathematical convenience is to convert
$\mathbb{C}\cup\{\infty\}\setminus\overline{\mathbb{D}}$ into
$\mathbb D=\{z\in\mathbb{C}:|z|<1\}$ by noting that $p_n(z)=z^nq_n(1/z)$,
$z\not\in\overline{\mathbb{D}}$, where for all $z\in\mathbb{D}$,
\[
  q_n(z) = \det\left(1- z\frac{A_n}{\sqrt n}\right)
\]
is the reciprocal polynomial of the characteristic polynomial $p_n$. Let
$\mathrm{H}(\mathbb{D})$ be the set of holomorphic or complex analytic
functions on $\mathbb{D}$, equipped with the topology of uniform convergence
on compact subsets, \emph{the compact-open topology}, see for instance
\cite{Shirai}. This allows to see $q_n$ as a random variable on
$\mathrm{H}(\mathbb{D})$ and gives a meaning to convergence in law of $q_n$ as
$n\to\infty$, namely, 
$q_n$ converges in law to some
random element $q$ of $\mathrm{H}(\mathbb{D})$
if for every bounded real
continuous function 
$f$ on $\mathrm{H}(\mathbb{D})$,
$\mathbb E[f(q_n)] \to \mathbb E[f(q)]$.
\begin{theorem}[Convergence of reciprocal characteristic polynomial]\label{th:q}
  We have
  \[
    q_n
    \xrightarrow[n \to \infty]
    {\mathrm{law}}\kappa \mathrm{e}^{-F},
  \]
  where $F$ is the random holomorphic function on $\mathbb D$ defined by
  \[
    F(z)=\sum_{k=1}^\infty X_k \frac{z^k}{\sqrt k}
  \]
  where $\{X_k\}_{k\geq1}$ is a sequence of independent complex Gaussian
  random variables such that
  \[
    \mathbb E\big[X_k\big]=0,
    \quad
    \mathbb E\left[|X_k|^2\right] = 1
    \quad\text{and}\quad
    \mathbb E\left[X_k^2 \right] = \mathbb E \left[a_{11}^2 \right]^k,
  \]
  and where $\kappa: \mathbb D \to \mathbb C$ is the holomorphic
  function %such that $\kappa(0)=1$ and
  defined for all $z\in\mathbb{D}$ by
  \[
    \kappa(z) = \sqrt{1-z^2 \mathbb E \left[a_{11}^2 \right]}.
  \]
\end{theorem}

The square root defining $\kappa$ is the one
such that $\kappa(0)=1$.
Notice that it is a well-defined holomorphic function on
the simply connected domain $\mathbb{D}$
since the function
$z\mapsto1-z^2\mathbb{E}[a_{11}^2]$ does not vanish on $\mathbb{D}$ which is true due to the
fact that $|\mathbb{E}[a_{11}^2]|\leq\mathbb{E}[|a_{11}|^2]=1$.

The proof of Theorem \ref{th:q} is given in Section \ref{se:pr:th:q}. It is
partially inspired by \cite{basak2019outliers} and relies crucially on a
joint combinatorial central limit theorem for traces of fixed powers (Lemma
\ref{le:clt}) inspired from \cite{JansonNowicki}. Unlike previous arguments
used in the literature for the analysis of Girko matrices, the approach does
not rely on Girko Hermitization, Gelfand spectral radius formula, high order
traces, resolvent method or Cauchy\,--\,Stieltjes transform. The first step
consists in showing the tightness of ${(q_n)}_{n\geq1}$, by using a
decomposition of the determinant into orthogonal elements related to
determinants of submatrices, as in \cite{basak2019outliers}. Knowing this
tightness, the problem is reduced to show the convergence in law of these
elements. A reduction step, inspired by \cite{JansonNowicki}, consists in
truncating the entries, reducing the analysis to the case of bounded entries.
The next step consists in a central limit theorem for product of traces of
powers of fixed order. It is important to note that we truncate with a fixed
threshold with respect to $n$, and the order of the powers in the traces are
fixed with respect to $n$. This is in sharp contrast with the usual
Füredi\,--\,Komlós truncation-trace approach related to the Gelfand spectral
radius formula used in \cite{MR866352,MR841088,MR863545,MR3813992}.

\subsection{Comments and open problems}

\subsubsection{Moment assumptions.}

The universality for the first order global asymptotics \eqref{eq:cirlaw}
depends only on the trace $\mathbb{E}[|a_{11}|^2]$ of the covariance matrix of
$\Re a_{11}$ and $\Im a_{11}$. The universality stated by Theorem \ref{th:q},
just like for the central limit theorem, depends on %more aspects of this
the whole covariance matrix. Since
\[
   \mathbb{E}[a_{11}^2]=\mathbb{E}[(\Re a_{11})^2]-\mathbb{E}[(\Im a_{11})^2]
  +2\mathrm{i}\mathbb{E}[\Re a_{11}\Im a_{11}],
\]
we can see that $\mathbb{E}[a_{11}^2]=0$ if and only if
$\mathbb{E}[(\Re a_{11})^2]=\mathbb{E}[(\Im a_{11})^2]$ and
$\mathbb{E}[\Re a_{11}\Im a_{11}]=0$. %\\
Moreover, we cannot in general get rid of $\mathbb{E}[a_{11}^2]$ by simply
multiplying the matrix $A_n$ by a phase.

\subsubsection{Hyperbolic Gaussian analytic function}

When $\mathbb E\left[a_{11}^2\right]=0$ then $\kappa=1$ while the random
analytic function $F$ which appears in the limit in Theorem \ref{se:pr:th:q}
is a degenerate case of the well-known hyperbolic Gaussian Analytic Functions
(GAFs) \cite[Equation (2.3.5)]{MR2552864}. It can also be obtained as the
antiderivative of the $L=2$ hyperbolic GAF which is $0$ at $z=0$. This $L=2$
hyperbolic GAF is related to the Bergman kernel and could be called the
Bergman GAF. These GAFs appear also at various places in mathematics and
physics and, in particular, in the asymptotic analysis of Haar unitary matrices,
see \cite{MR1844632,hmc}.

\subsubsection{Cauchy\,--\,Stieltjes transform}

If $\mathbb{E}[a_{11}^2]=0$ then by returning to $p_n$, taking the logarithm
and the derivative with respect to $z$ in Theorem \ref{th:q}, we obtain the
convergence in law of the Cauchy\,--\,Stieltjes transform (complex conjugate
of the electric field) minus $n/z$ towards $z \mapsto F'(1/z)/z^2$ which is a
Gaussian analytic function on $\mathbb C \setminus \overline{\mathbb D}$ with
covariance given by a Bergman kernel.

\subsubsection{Central Limit Theorem}
We should see Theorem~\ref{th:q} as a global second order analysis, just like
the central limit theorem (CLT) for linear spectral statistics
\cite{MR2294978,cipolloni2020fluctuation,cipolloni2019central}.
Namely for all $z\in\mathbb{D}$, we have 
$|q_n(1/z)| = \exp
\left[-n \left(U_n(z)-U(z)\right)\right]$ where $U_n(z)=-\frac{1}{n}\log|p_n(z)|$ is the logarithmic potential at the 
point $z$ of the empirical spectral distribution of $\frac{1}{\sqrt n} A_n$ and $U(z)=-\log|z|$ is the logarithmic potential at the point $z$ of the uniform distribution on the unit disc $\mathbb{D}$.

Moreover, it is possible
to extract from Theorem~\ref{th:q} a CLT for linear spectral statistics with
respect to analytic functions in a neighborhood of $\overline{\mathbb{D}}$.
This can be done by using the Cauchy formula for an analytic function $f$,
\[
  \int f(\lambda) \mathrm \mu(\mathrm d \lambda)
  = \frac{1}{2\pi\mathrm{i}}
    \int \left(\oint \frac{f(z)}{z - \lambda} \mathrm d z\right) \mu(\mathrm d
    \lambda)
  = \frac{1}{2\pi\mathrm{i}} \oint f(z) \left(\int \frac{\mu(\mathrm d\lambda)} {z - \lambda} \right)\mathrm d z\\%
  = \frac{1}{2\pi\mathrm{i}} \oint f(z)(\log \det \left(z - A \right))' \mathrm d z%
\]
where $\mu$ is the counting measure of the eigenvalues of $A$, where the
contour integral is around a centered circle of radius strictly larger than
$1$, and where we have taken any branch of the logarithm. The approach is
purely complex analytic. In particular, it is different from the usual
approach with the logarithmic potential of $\mu$ based on the real function
given by
$z \mapsto \int\log|z-\lambda|\mu(\mathrm{d}\lambda) = \log|\det(z-A)|$.
% It involves the complex Cauchy\,--\,Stieltjes transform
% $z\mapsto\int\frac{\mu(\mathrm{d}\lambda)}{z-\lambda}$ of $\mu$. , which is
% locally Lebesgue integrable.

\subsubsection{Wigner case and elliptic interpolation}

The finite second moment assumption of Theorem \ref{th:p} is optimal. We could
explore its relation with the finite fourth moment assumption for the
convergence of the spectral edge of Wigner random matrices, which is also
optimal. Heuristic arguments tell us that the interpolating condition on the matrix entries 
should be
$\mathbb{E}[|a_{jk}|^2|a_{kj}|^2]<\infty$ for $j\neq k$, which is a finite
second moment condition for Girko matrices and a finite fourth moment
condition for Wigner matrices. 
This is work in progress.

\subsubsection{Coupling and almost sure convergence}

For simplicity, we define in \eqref{eq:An} our random matrix $A_n$ for all
$n\geq1$ by truncating from the upper left corner the infinite random matrix
$\{a_{jk}\}_{j,k\geq1}$. This imposes a coupling for the matrices
$\{A_n\}_{n\geq1}$. However, since Theorem \ref{th:p} involves a convergence
in probability, it remains valid for an arbitrary coupling, in the spirit of
the triangular arrays assumptions used for classical central limit theorems.
In another direction, one could ask about the upgrade of the convergence in
probability into almost sure convergence in Theorem \ref{th:p}.
This is an open problem.

\subsubsection{Heavy tails}

An analogue of \eqref{eq:cirlaw} in the heavy-tailed case
$\mathbb{E}[|a_{11}|^2]=\infty$ is considered in \cite{MR2837123} but requires
another scaling than $\frac{1}{\sqrt{n}}$. The spectral radius of this model
tends to infinity as $n\to\infty$ but it could be possible to analyze the
limiting point process at the edge as $n\to\infty$ and its universality.
This is an open problem.

\section{Proof of Theorem \ref{th:p}}
\label{se:pr:th:p}

Let $f=\kappa\mathrm{e}^{-F}$ be as in Theorem \ref{th:q}. We observe that the
equation $f(z)=0$, $z\in\mathbb{D}$ is equivalent to $\kappa(z)=0$,
$z\in\mathbb{D}$, which has no solution, because
$|\mathbb{E}[a_{11}^2]|\leq\mathbb{E}[|a_{11}|^2]=1$. In particular, for every
$r\in(0,1)$,
\[
  \inf_{z\in\overline D_r}|f(z)|
  =\inf_{z \in \overline D_r} \left\{|\kappa(z)|\mathrm{e}^{-\Re(F(z))}\right\}>0  ,
\]
where $\overline D_r=\{z\in\mathbb{C}:|z|\leq r\}$ is the closed disc of
radius $r$. On the other hand, the convergence in law provided by Theorem
\ref{th:q} together with the continuous mapping theorem 
used for
the continuous function $f \in \mathrm{H}(\mathbb{D})
\mapsto \inf_{z \in \overline D_r} |f(z)|$ give, for every
$r\in(0,1)$,
\[
  \inf_{z \in \overline D_r}
  |q_n(z)|
  \xrightarrow[n \to \infty]{\mathrm{law}} 
  \inf_{z \in \overline D_r}
  \left\{|\kappa(z)|\mathrm{e}^{-\Re(F(z))}
  \right\}.
\]
Now, since $q_n(z)=z^np_n(1/z)$ for every $z\in\mathbb{D}$,
we obtain, by combining these two facts, for every $r\in(0,1)$,
\[
  \mathbb P\Bigr(\rho_n<\frac{1}{r}\Bigr)
  = \mathbb P\Bigr(\inf_{|z|\geq\frac{1}{r}}|p_n(z)|>0\Bigr)
  = \mathbb P\Bigr(\inf_{z \in \overline D_r}|q_n(z)|>0\Bigr)
  \xrightarrow[n\to\infty]{}
  \mathbb P\Bigr(\inf_{z \in \overline D_r}\left\{|\kappa(z)|\mathrm{e}^{-\Re(F(z))}\right\}>0\Bigr)=1.
\]
In other words, for all $\varepsilon>0$,
\[
  \lim_{n\to\infty}\mathbb{P}(\rho_n\geq1+\varepsilon)
  =0.
\]
Combined with \eqref{eq:liminf}, this leads to the desired result
\[
  \lim_{n\to\infty}\mathbb{P}(|\rho_n-1|>\varepsilon)
  =0.
\]
%Note that there is also almost sure convergence for a coupling provided by Skorokhod representation.

Note that it could be possible to obtain the result by using the Radon
measures of the zeros and the Hurwitz phenomenon, see \cite[Lemma 2.2]{Shirai}
and \cite[Lemma 5.2]{raphael-david}, but this would be more complicated!

\section{Proof of Theorem \ref{th:q}}
\label{se:pr:th:q}

By developing the determinant we can see that
\[
  q_n(z)
  =\det\left(1-
    z\frac{A_n}{\sqrt n}\right)
  = 
  1 + \sum_{k=1}^n (-z)^{k} P^{(n)}_k,
\]
where
\[
  P^{(n)}_k = \sum_{\substack{I \subset \{1,\dots,n\}\\|I|=k}}n^{-k/2} \det(A_n(I))
  \quad\text{and}\quad
  A_n(I) = \{a_{jk}\}_{j,k \in I}.
\]

The following lemma is essentially contained in \cite[Appendix A]{basak2019outliers}. It
 is proved in Section \ref{se:le:tightness}.

\begin{lemma}[Tightness]\label{le:tightness}
  The sequence $\{q_n\}_{n\geq1}$ is tight.
\end{lemma}

For completeness, let us recall that 
the sequence $\{q_n\}_{n\geq1}$
 is tight if for every $\varepsilon> 0$,
 there exists a compact subset of $\mathrm{H}(\mathbb D)$
 such that $\mathbb P(q_n \in K)> 1-\varepsilon$
 for every $n$.

Now that we know that $\{q_n\}_{n\geq1}$ is tight, it is enough to understand,
for each $k \geq 1$, the limit of $(P_1^{(n)},\dots,P_k^{(n)})$ as
$n \to \infty$. Indeed, we have the following Lemma \ref{le:coefcv}, close to
\cite[Second part of Proposition 2.5]{Shirai}. For the reader's convenience and
for completeness, we give a proof in Section \ref{se:le:coefcv}.

\begin{lemma}[Reduction to convergence of coefficients]\label{le:coefcv}
  Let $\{f_n\}_{n\geq1}$ be a tight sequence of random elements of
   $\mathrm{H}(\mathbb{D})$, and let us write, for every $n\geq1$,
  $f_n(z)=\sum_{k=0}^\infty z^k P_k^{(n)}$. If for every $m \geq 0$,
  \[
    (P_0^{(n)},\dots,P_m^{(n)})
    \xrightarrow[n \to \infty]{\mathrm {law}}
    (P_0,\dots,P_m)
  \]
  for a 
  common 
  sequence of random variables $\{P_m\}_{m\geq 0}$ then
  $f=\sum_{k=0}^\infty z^kP_k$ is well-defined in $\mathrm{H}(\mathbb{D})$ and
  \[
    f_n
    \xrightarrow[n \to \infty]{\mathrm {law}}
    f.
  \]
\end{lemma}

The first simplification we shall make is to assume that $a_{11}$ is bounded.
This is motivated by \cite[Proof of Lemma 7]{JansonNowicki}. We write this in
the following lemma, proved in Section \ref{se:le:truncation}.

\begin{lemma}[Reduction to bounded entries by truncation]\label{le:truncation}
  For $M>0$ let us define
  \[
    A_n^{(M)}= \left\{a_{ij}^{(M)}\right\}_{1 \leq i,j \leq n}
    \quad\text{where}\quad    
    a_{ij}^{(M)} %
    = a_{ij} 1_{|a_{ij}| < M} %
    - \mathbb E \left[a_{ij} 1_{|a_{ij}| < M}\right]
  \]
  and 
  \[
    P^{(n,M)}_k = \sum_{\substack{I \subset \{1,\dots,n\}\\|I|=k}} n^{-k/2} \det(A_n^{(M)}(I))
    \quad\text{where}\quad
    A_n^{(M)}(I) = \left\{a_{ij}^{(M)}\right\}_{i,j \in I}.
  \]
  Let $k \geq 1$. If there exists 
  $\left\{\left(Y_1^{(M)},\dots,Y_k^{(M)} \right) \right\}_{M\geq1}$ and a
  random vector $\left(Y_1,\dots,Y_k \right)$ such that for all $M\geq1$,
  \[
    \left(P_1^{(n,M)}, \dots,P_k^{(n,M)} \right)\xrightarrow[n \to
    \infty]{\mathrm{law}} \left(Y_1^{(M)},\dots,Y_k^{(M)} \right),
    \quad\text{and}\quad
    \left(Y_1^{(M)},\dots,Y_k^{(M)} \right) \xrightarrow[M \to
    \infty]{\mathrm{law}} \Big(Y_1,\dots,Y_k \Big),
  \]
  then 
  \[
    (P_1^{(n)},\dots,P_k^{(n)}) 
    \xrightarrow[n \to \infty]{\mathrm{law}} 
    (Y_1,\dots,Y_k).
  \]
\end{lemma}

To simplify the study of $P_k^{(n)}$ we notice the following. For each integer 
$n \geq 1$, the series
\[
  B_n = -\sum_{k=1}^\infty 
  \left(\frac{A_n}{\sqrt n} \right)^k
  \frac{z^k}{k}
\]
converges for $z$ small enough and its exponential is
$\left(1- z\frac{A_n}{\sqrt n} \right)$.
This can be shown in the standard way
if $A_n$ is diagonalizable and can be extended
to non-diagonalizable matrices by continuity. Then, since
$\det\left(\mathrm{e}^{B_n}\right)=\mathrm{e}^{\mathrm{Tr}B_n}$, we obtain
\[
  q_n(z) = \exp\left(-\sum_{k=1}^\infty\frac{\mathrm{Tr}(A_n^k)}{n^{k/2}}
    \frac{ z^k}{k}
  \right)
\]
for $z$ small enough.
In particular, $(P_1^{(n)},\dots,P_k^{(n)})$ is a polynomial function of
$\left( \frac{\mathrm{Tr}(A_n)}{n^{1/2}}, \dots,\frac{\mathrm{Tr}(A_n^k)}
  {n^{k/2}} \right)$ that does not depend on $n$ and vice versa. The idea is
to study, by the method of moments, the quantity
\[
  \mathrm{Tr}(A_n^k)
  = \sum_{
    1 \leq i_1,\dots,i_k \leq n}
  a_{i_1 i_2}
  a_{i_2 i_3}\dots a_{i_{k-1} i_k}
  a_{i_k i_1}.
\]
That is why we preferred to have $a_{11}$ bounded (or at least having all
its moments finite). Note that we have used the determinantal terms $P_k^{(n)}$ to perform this truncation step, it would have been much more challenging to justify this truncation directly for the traces $\mathrm{Tr}(A_n^k)$. On the other hand, it would have been much more difficult to prove directly the convergence of the determinantal terms $P_k^{(n)}$ thanks to the method of moments since these terms are  asymptotically neither independent nor Gaussian.

We decompose the above sum in two sums,
\begin{equation}\label{eq:TrSum}
  \mathrm{Tr}(A_n^k)
  = 
  \sum_{\substack{1 \leq i_1,\dots,i_k \leq n\\\mathrm{card}\left\{i_1,\dots,i_k\right\} = k}}
  a_{i_1 i_2}
  a_{i_2 i_3}\dots a_{i_{k-1} i_k}
  a_{i_k i_1}
  +
  \sum_{\substack{1 \leq i_1,\dots,i_k \leq n\\\mathrm{card} \left\{i_1,\dots,i_k\right\} < k}}
  a_{i_1 i_2}
  a_{i_2 i_3}\dots a_{i_{k-1} i_k}
  a_{i_k i_1}.
\end{equation}
The first term in the right-hand side of \eqref{eq:TrSum} has zero expected
value and gives rise to the random part of the limit. The second term in the
right-hand side of \eqref{eq:TrSum} gives the deterministic part.

We begin by looking at the term
\begin{equation}\label{eq:TrDifferent}
  \sum_{\substack{1 \leq i_1,\dots,i_k \leq n\\\mathrm{card} \left\{i_1,\dots,i_k\right\} = k}}
  a_{i_1 i_2}
  a_{i_2 i_3}\dots a_{i_{k-1} i_k}
  a_{i_k i_1}.
\end{equation}
Notice that the sum in \eqref{eq:TrDifferent} is indexed by sequences
$\iota = (i_1,\dots,i_k)$ of pairwise distinct elements of $\{1,\dots,n\}$.
However, if two sequences are cyclic permutations of each other we obtain the
same term. To deal with this fact, we should consider sequences up to cyclic
permutations or, what is the same, directed cycles in $\{1,\dots,n\}$. More
precisely, let us consider $\{1,\dots,n\}$ as the complete directed graph with
no loops and let us consider the graph $G$ = $(V,E)$ with vertex and edge sets
\[
  V=
  \{1,\dots,k\} \quad \text{and} \quad 
  E = \{(1,2),(2,3),\dots,(k-1,k),(k,1)\}.
\]
A $k$-directed cycle in $\{1,\dots,n\}$ is a subgraph $g$ of $\{1,\dots,n\}$
that is isomorphic to $G$. The sum in \eqref{eq:TrDifferent} is better indexed
by the set of $k$-directed cycles in $\{1,\dots,n\}$ that we shall call
$\mathcal C_k^{(n)}$. For $g \in \mathcal C_k^{(n)}$, we define
\[
  a_g = \prod_{e\text{ edge of }g} a_e,
\]
where $a_e = a_{ij}$ if $e = (i,j)$. Now, we can write
\[
  \sum_{\substack{1 \leq i_1,\dots,i_k \leq n\\\mathrm{card} \left\{i_1,\dots,i_k\right\} = k}}
  a_{i_1 i_2}
  a_{i_2 i_3}\dots a_{i_{k-1} i_k}
  a_{i_k i_1}
  = k
  \sum_{g \in \mathcal C_k^{(n)}}
  a_{g}
\]
so that the term we have to study is
\[
  t_k^{(n)}
  = \sum_{g \in \mathcal C_k^{(n)}}
  a_{g}.
\]

The following lemma is a sort of combinatorial joint central limit theorem. It
provides the $\mathrm{e}^{-F}$ part of the limiting random analytic function
in Theorem \ref{th:q}. It is proved in Section \ref{se:le:clt}.

\begin{lemma}[Convergence to a Gaussian object]\label{le:clt}
  For any $k_1,\dots,k_m \geq 1$ and any sequence
  $s_1,\dots,s_m \in \{\cdot,*\} $,
  \[
    \mathbb E \left[
      \left(\frac{
          t_{k_1}^{(n)}}{n^{k_1/2}}\right)^{s_1}
      \cdots
      \left(\frac{t_{k_m}^{(n)}}
        {n^{k_m/2}}\right)^{s_m}
    \right]
    \xrightarrow[n \to \infty]{}
    \mathbb E
    \left[
      \left(\frac{X_{k_1}}{\sqrt{k_1}
        }\right)^{s_1}\cdots 
      \left(\frac{X_{k_m}}
        {\sqrt{k_m}} \right)^{s_m} \right],
  \]
  where we have used the notation $x^\cdot = x$ and $x^* = \bar x$, and where
  $\{X_k\}_{k\geq1}$ are independent complex Gaussian random variables such
  that $\mathbb E\big[X_k\big]=0$, $\mathbb E\left[|X_k|^2\right] = 1$, and
  $\mathbb E\left[X_k^2 \right] = \mathbb E \left[a_{11}^2 \right]^k$ for all
  $k\geq1$.
\end{lemma}

The term that is left to understand is
\[
  r_k^{(n)}=
  \sum_{\substack{1 \leq i_1,\dots,i_k \leq n\\\mathrm{card} \left\{i_1,\dots,i_k\right\} < k}}
  a_{i_1 i_2}
  a_{i_2 i_3}\cdots a_{i_{k-1} i_k}
  a_{i_k i_1}.
\]

The following lemma, proved in Section \ref{se:le:meancv}, provides the
$\kappa$ part of the limit in Theorem \ref{th:q}.

\begin{lemma}[Deterministic limit part]\label{le:meancv}
  \[
    \frac{r_k^{(n)}}{n^{k/2}}
    \xrightarrow[n \to \infty]
    {\mathrm{law}} 
    \begin{cases}
      \mathbb E\left[{a_{11}}^2\right]^{k/2} & \text{if $k$ is even}\\
      0 & \text{if $k$ is odd.}
    \end{cases}
  \]
\end{lemma}
To sum up, if $\{X_k\}_{k\geq1}$ is a sequence of independent complex Gaussian
random variables such that
\[
  \mathbb E\left[|X_k|^2\right] = 1
  \quad\text{and}\quad
  \mathbb E\left[X_k^2 \right] = \mathbb E \left[a_{11}^2 \right]^k,
\]
and if
\[
  \text{mean}_k
  =\begin{cases}
    \mathbb E\left[a_{11}^2 \right]^{k/2}&\text{if $k$ is even}\\
    0&\text{if $k$ is odd}
  \end{cases}
\]
then
\[
  \left(\frac{\mathrm{Tr}(A_n)}{\sqrt n},
    \frac{\mathrm{Tr}(A_n^2)}{n},
    \dots,
    \frac{
      \mathrm{Tr}(A_n^k)}{n^{k/2}}
  \right)
  \xrightarrow[n \to \infty]{\mathrm{law}}
  \left(X_1,\sqrt 2 X_2,
    \dots, \sqrt k X_k  \right)
  +
  \left(\text{mean}_1,\text{mean}_2,\dots,\text{mean}_k \right).
\]
This implies the convergence of $(P_1^{(n)},\dots,P_k^{(n)})$ to the
corresponding polynomials of $X_i$ and $\text{mean}_i$. Moreover, by Lemma
\ref{le:coefcv} and since the limit depends continuously on the second moment
of the variable, the assertion also holds for non-bounded $a_{11}$. We have
found that $\{q_n\}_{n\geq1}$ has a limit that can be written as a Maclaurin
series whose coefficients are polynomials of $X_i$ and $\text{mean}_i$. By
construction, the joint law of these coefficients is the same as the joint law
of the coefficients of the random holomorphic function $\kappa \exp(-F)$ so
that the proof of the theorem is complete.

\section{Proofs of the lemmas used in the proof of Theorem \ref{th:q}}
\label{se:lemmasproofs}

\subsection{Proof of Lemma \ref{le:tightness}}\label{se:le:tightness}

Recall that if $\{f_n\}_{n\geq1}$ is a sequence of random variables on
$\mathrm{H}(\mathbb{D})$ such that for every compact set $K\subset\mathbb{D}$
the sequence of random variables $\{\Vert f_n\Vert_K\}_{n\geq1}$ is tight,
$\Vert f_n\Vert_K=\max_K|f_n|$, then $\{f_n\}_{n\geq1}$ is tight, see for
instance \cite[Proposition 2.5]{Shirai}.

By \cite[Lemma 2.6]{Shirai}, it is enough to bound $\mathbb E[|q_n(z)|^2]$ by a deterministic continuous
function of $z$ that does not depend on $n$. 
Recall the notation $A_n(I)$ from the beginning of Section \ref{se:pr:th:q}. Notice that each $\det(A_n(I))$
has mean zero and that
\[
  \mathbb E \left[\det(A_n(I)) \overline{\det( A_n(J))} \right]=0
  \quad\text{if } I \neq J
  \quad\text{and}\quad
  \mathbb E[\left|\det(A_n(I)) \right|^2]= \mathrm{card}(I)!.
\]
In particular, 
\[\mathbb E\left[|P^{(n)}_k|^2 \right]
  = n^{-k}
  \binom{n}{k}
  \mathbb E\left[|\det A_k|^2 \right]
  =n^{-k} \frac{n!}{(n-k)!} \leq 1 
  \quad\text{and}\quad
  \mathbb E
  \left[P_k^{(n)}\overline{P_l^{(n)}} \right]=0 
  \quad\text{if}\quad k \neq l.
\]
So, we have
\[
  \mathbb E[|q_n(z)|^2]
  \leq 1 + \sum_{k=1}^n |z|^{2k} \mathbb E[|P^{(n)}_k|^2]
  \leq \sum_{k=0}^n |z|^{2k}
  \leq \frac{1}{1-|z|^2}.
\]

\subsection{Proof of Lemma \ref{le:coefcv}}
\label{se:le:coefcv}

The statement is close to \cite[Proposition 2.5]{Shirai}.

Take two subsequences $\{f_{n_\ell}\}_{\ell\geq1}$ and
$\{f_{\tilde n_\ell}\}_{\ell\geq1}$ of random functions that converge, in law,
to some random functions $g$ and $\widetilde g$ in $\mathrm{H}(\mathbb{D})$.
We want to show that the distributions of $g$ and $\widetilde g$ coincide. By
Remark \ref{rk:zeros} below, we can write $g(z)=\sum_{k=0}^\infty Q_kz^k$ for
$z\in\mathbb{D}$ and $\widetilde g(z)= \sum_{k=0}^\infty \widetilde Q_kz^k$
for $z\in\mathbb{D}$, where $\{Q_k\}_{k\geq0}$ and
$\{\widetilde Q_k\}_{k\geq0}$ are two sequences of complex random variables.
By the same remark, we have that for any $m \geq 0$, the limit in law as
$\ell\to\infty$ of $(P_0^{(n_\ell)},\dots,P_m^{(n_\ell)})$ is
$(Q_0,\dots,Q_m)$ while the limit in law as $\ell\to\infty$ of
$(P_0^{(\tilde n_\ell)},\dots,P_m^{(\tilde n_\ell)})$ is
$(\widetilde Q_0,\dots,\widetilde Q_m)$. In particular,
$(\widetilde Q_0,\dots, \widetilde Q_m)$ and $(Q_0,\dots,Q_m)$ have the same
distribution as $(P_0,\dots,P_m)$ for every $m \geq 0$ so that
$\{Q_k\}_{k \geq 0}$, $\{\widetilde Q_k\}_{k \geq 0}$ and $\{P_k\}_{k \geq 0}$
have the same distribution as random elements of
$\mathbb C^{\mathbb Z_{\geq 0}}$. By Remark \ref{rk:zeros} again, $g$ and
$\widetilde g$ have also the same distribution. Moreover, the random function
$z\in\mathbb{D}\mapsto f(z)=\sum_{k\geq 0}P_kz^k$ is well-defined as a random
variable on $\mathrm{H}(\mathbb{D})$ and its distribution is the unique limit
point of the sequence of distributions of $\{f_n\}_{n \geq 1}$. Finally, since
$\{f_n\}_{n\geq1}$ is tight and since, by Prokhorov's theorem, tightness means
that its sequence of distributions is sequentially relatively compact in the space of probability measures on 
$\mathrm{H}(\mathbb{D})$, we conclude that $\{f_n\}_{n\geq1}$
converges in law to $f$ as $n\to\infty$.

\begin{remark}\label{rk:zeros}
  The $P_k$'s are related to the successive derivatives of $f$ at point $0$.
  Due to the properties of analytic functions, the map
  $T:H(\mathbb D) \to \mathbb C^{\mathbb Z_{\geq 0}}$ defined for all
  $h\in\mathrm{H}(\mathbb{D})$ and all $k\in\mathbb{Z}_{\geq0}$ by
  \[
    T(h)_{k} = \frac{1}{k!}\frac{\mathrm d^k h}{\mathrm d z^k}(0)
  \]
  is continuous and injective. The inverse map
  $T^{-1}:\bigr\{\{a_k\}_{k\geq0} \in \mathbb C^{\mathbb Z_{\geq 0}}:
  \varlimsup_{k\to\infty}|a_k|^{1/k}\leq1\bigr\} \to H(\mathbb D)$ given by
  \[
    (T^{-1}(a))(z)=\sum_{k=0}^\infty a_kz^k
  \]
  is measurable. Denoting $\mathcal{P}(E)$ the set of probability measures on
  $E$, it follows that the pushforward map
  \[
    T_*: \mathcal P\big(H(\mathbb D)\big) \to \mathcal P\left(\mathbb C^{\mathbb
        Z_{\geq 0}} \right)
  \]
  is injective in the sense that for all $\mu$ and $\nu$ in
  $\mathcal{P}(\mathrm{H}(\mathbb{D}))$, if $T_*\mu = T_*\nu$ then $\mu =
  \nu$.
\end{remark}

\subsection{Proof of Lemma \ref{le:truncation}}\label{se:le:truncation}

It is enough to notice that, for each $k \geq 1$, there exists a sequence
$\left\{C_M \right\}_{M\geq1} $ that goes to zero such that
\[
  \mathbb E\Bigr[\Bigr|P_k^{(n,M)} - P_k^{(n)}\Bigr|^2\Bigr] \leq C_M
\]
for every $n, M \geq 1$. But
\begin{align*}
  \mathbb E\Bigr[\Bigr|P_k^{(n,M)} - P_k^{(n)}\Bigr|^2\Bigr]
  &=
    n^{-k}
    \sum_{\substack{I \subset \{1,\dots,n\}\\|I|=k}}
    \mathbb E \left[
    \left|\det(A_n(I)^{(M)}) -
    \det(A_n(I))\right|^2 \right]					\\
  &=
	n^{-k}\binom{n}{k} 
    \mathbb E\Bigr[\Bigr|a_{11}^{(M)}\cdots a_{1k}^{(M)} - a_{11}\cdots a_{1k}\Bigr|^2\Bigr]		
    k!												\\
  &\leq
    \mathbb E\Bigr[\Bigr|a_{11}^{(M)}\cdots a_{1k}^{(M)} - a_{11}\cdots a_{1k} \Bigr|^2\Bigr]		
\end{align*}
so that
$C_M = \mathbb E\Bigr[\Bigr|a_{11}^{(M)}\cdots a_{1k}^{(M)} - a_{11}\cdots a_{1k}\Bigr|^2\Bigr]$ works.

\subsection{Proof of Lemma \ref{le:clt}}
\label{se:le:clt}

As it is usual, the idea is to understand which terms are dominant. We have
\begin{align*}
  \mathbb E \left[
  \left(\frac{t_{k_1}^{(n)}}{n^{k_1/2}}\right)^{s_1}\cdots\left(\frac{t_{k_m}^{(n)}}{n^{k_m/2}}\right)^{s_m}
  \right]
  &=
    \frac{1}{n^{(k_1+\dots+k_m)/2}}
    \mathbb E \left[
    \left(t_{k_1}^{(n)}\right)^{s_1}
    \cdots
    \left(t_{k_m}^{(n)}\right)^{s_m}
    \right]							\\
  &=
    \frac{1}{n^{(k_1+\dots+k_m)/2}}
    \sum_{g_1\in\mathcal C_{k_1}^{(n)},\ldots, g_m \in \mathcal C_{k_m}^{(n)}}
    \mathbb E\left[\left(a_{g_1}\right)^{s_1}\cdots\left(a_{g_m}\right)^{s_m}
    \right].
\end{align*}
We say that $(g_1,\dots,g_m)$ is equivalent to
$(\widetilde g_1,\dots, \widetilde g_m)$ if there is a bijection
$\theta:\{1,\dots,n\} \to \{1,\dots,n\}$ such that
\[
  \widetilde g_i = \theta_* (g_i)
  \quad \text{for } i \in 
  \{1,\dots,m\},
\]
where $\theta_{*}$ denotes the map induced by $\theta$ on the subgraphs
of $\{1,\dots,n\}$. So,
\[
  \mathbb E\left[
    \left(a_{ g_1}
    \right)^{s_1} \cdots 
    \left(a_{ g_m}\right)^{s_m}
  \right]
  =
  \mathbb E\left[
    \left(a_{\widetilde g_1}
    \right)^{s_1} \cdots 
    \left(a_{\widetilde g_m}\right)^{s_m}
  \right]
\]
if $\Gamma= (g_1,\dots,g_m)$ is equivalent to
$\widetilde \Gamma = (\widetilde g_1,\dots, 
\widetilde g_m)$. Hence, if we denote by
$\mathcal T_{(k_1,\dots,k_m)}^{(n)}$ the set of equivalence classes, we can
define
\[
  W_{[\Gamma]}=\mathbb E\left[ \left(a_{g_1} \right)^{s_1} \cdots \left(a_{g_m}\right)^{s_m} \right],
\]
where $[\Gamma]$ is the class of $\Gamma$. We can then write
\[
  \frac{1}{n^{(k_1+\dots+k_m)/2}}
  \sum_{g_1 \in \mathcal C_{k_1}^{(n)},\ldots, g_m \in \mathcal C_{k_m}^{(n)}}
  \mathbb E\left[
    \left(a_{g_1}
    \right)^{s_1} \cdots 
    \left(a_{g_m}\right)^{s_m}
  \right]
  =
  \frac{1}{n^{(k_1+\dots+k_m)/2}}
  \sum_{\nu\in\mathcal T_{(k_1,\dots,k_m)}^{(n)}}
  \mathrm{card}(\nu)
  W_\nu,
\]
\begin{figure}[htbp]\centering
  \includegraphics[width=0.6\columnwidth]{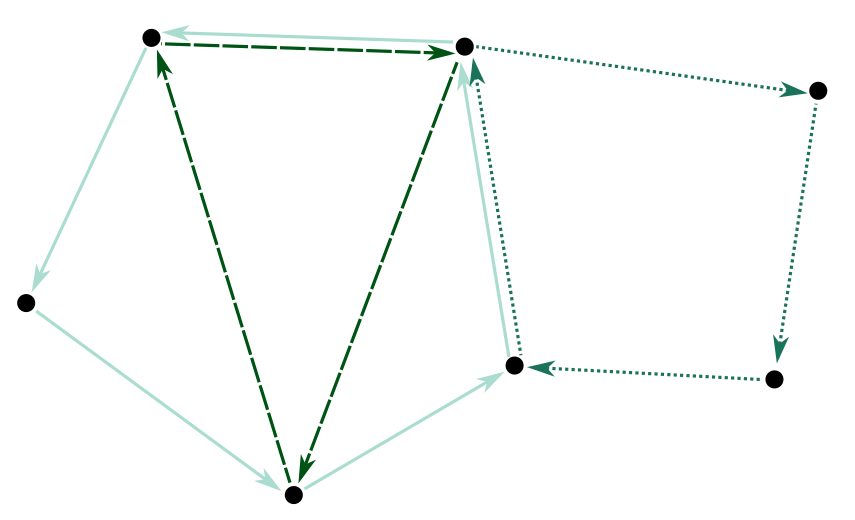}
  \caption{An example of a multigraph $E^\mu$ constructed from one
    $3$-directed cycle, one $4$-directed cycle and one $5$-directed
    cycle.\label{fig:GraphOfCycles}}
\end{figure}

\noindent
where $\mathrm{card}(\nu)$ is the cardinality of $\nu$ seen as a subset of
$\mathcal C_{k_1}^{(n)} \times \dots \times \mathcal C_{k_m}^{(n)}$. There is
a natural inclusion map from $\mathcal T_{(k_1,\dots,k_m)}^{(n)}$ into
$\mathcal T_{(k_1,\dots,k_m)}^{(n+1)}$ induced by the inclusion
$\{1,\dots,n\} \subset \{1,\dots,n+1\}$ and these inclusions are surjective if
$n \geq k_1 + \dots + k_m$. With the help of these inclusions we can write,
for $n\geq k_1+\dots+k_m$,
\[
  \frac{1}{n^{(k_1+\dots+k_m)/2}}
  \sum_{\nu\in 
    \mathcal T_{(k_1,\dots,k_m)}^{(n)}
  }
  \mathrm{card}(\nu)
  W_\nu
  =
  \frac{1}{n^{(k_1+\dots+k_m)/2}}
  \sum_{\mu\in 
    \mathcal T_{(k_1,\dots,k_m)}^{(k_1+
      \dots + k_m)}
  } 
  \! \! \! \! \!
  \mathrm{card}_n\left(\mu 
%    \cap 
%    \mathcal C_{k_1}^{(n)}
%    \times \dots \times 
%´    \mathcal C_{k_m}^{(n)}
  \right)
  W_\mu,
\]
where $\mathrm{card}_n\left(\mu \right)$
denotes the cardinality of
$\mu$ when seen as a subset of 
$\mathcal C_{k_1}^{(n)} \times \dots \times 
\mathcal C_{k_m}^{(n)}$.
So, it is enough to find the limit, as $n \to \infty$, of
\[\frac{
    \mathrm{card}_n\left(\mu 
%      \cap 
%      {\mathcal C}_{k_1}^{(n)}
%      \times \dots \times 
%      {\mathcal C}_{k_m}^{(n)}
    \right)}
  {n^{(k_1+\dots+k_m)/2}}\] for any
$\mu \in \mathcal T_{(k_1,\dots,k_m)}^{(k_1+ \dots + k_m)}$. To understand
better this cardinality, to each
\[
  \mu = [(g_1,\dots,g_m)] \in \mathcal T_{(k_1,\dots,k_m)}^{(k_1+ \dots +
    k_m)}
\]
we associate the oriented multigraph $G^\mu$ consisting of the union of the
$g_l$'s with edges counted multiple times (see Figure
\ref{fig:GraphOfCycles}). More precisely, if $V^{g_l}$ and $E^{g_l}$ are the
vertex set and the edge set of $g_l$, then the vertex set $V^{\mu}$ and the
edge set $E^{\mu}$ of $G^{\mu}$ are
\[
  V^{\mu} = \bigcup_{l=1}^m V^{g_l}
  \quad\text{and}\quad
  E^{\mu} = \bigcup_{l=1}^m \left(\{l\} \times E^{g_l} \right)
\]
with the source and target maps, $s: E^\mu \to V^\mu$ and
$t: E^\mu \to V^\mu$, defined by
\[
  s \big(l,(i,j)\big)
  = i \quad \text{and} \quad t \big(l,(i,j)\big)
  = j.
\]
If there is an edge $e \in E^\mu $ that is not multiple, in other words such
that $(s,t)(e) \neq (s,t)(e')$ for every other edge $e' \neq e$, then
$W_\mu = 0$. So we consider only graphs where all edges are multiple. If for
each $v \in V^{\mu}$ the outer degree $\mathrm{deg}(v)$ is defined by
\[
  \mathrm{deg}(v)
  = \mathrm{card}\{e \in E^\mu:
  s(e) = v\}
\]
we have that $\mathrm{deg}(v) \geq 2$ for every $v \in V^\mu$. By using the
handshaking lemma, we have
\[
  \sum_{v \in V^\mu} \mathrm{deg}(v) = \mathrm{card}(E^\mu) = k_1 +
  \dots + k_m.
\]
We notice that if, moreover, $\mathrm{deg}(v_*) \geq 3$ for some
$v_* \in V^\mu$ then
\[
  k_1 + \dots + k_m
  =\sum_{v \in V^\mu \setminus
    \{v_*\}} 
  \mathrm{deg}(v)
  + \mathrm{deg}(v_*)
  \geq 
  2\left(\mathrm{card}(V^\mu)-1
  \right)
  +3
  = 2
  \mathrm{card}(V^\mu) + 1\]
so that
\[
  \mathrm{card}(V^\mu) < \frac{k_1+\dots + k_m}{2}.
\]
But
%\[
%  \mathrm{card}\left(\mu 
%    \cap 
%    \mathcal C_{k_1}^{(n)}
%    \times \dots \times 
%    \mathcal C_{k_m}^{(n)}
%  \right) \leq
$\mathrm{card}_n(\mu) \leq n^{\mathrm{card}(V^\mu)}$,
%\]
which implies that
\[
  \frac{\mathrm{card}_n(\mu)
%    \cap
%      \mathcal C_{k_1}^{(n)}
%      \times \dots \times 
%      \mathcal C_{k_m}^{(n)} 
}{n^{(k_1+\dots +k_m)/2}}
  \xrightarrow[n \to \infty]{} 0.
\]
\begin{figure}[htbp]\centering
  \includegraphics[width=0.6\columnwidth]{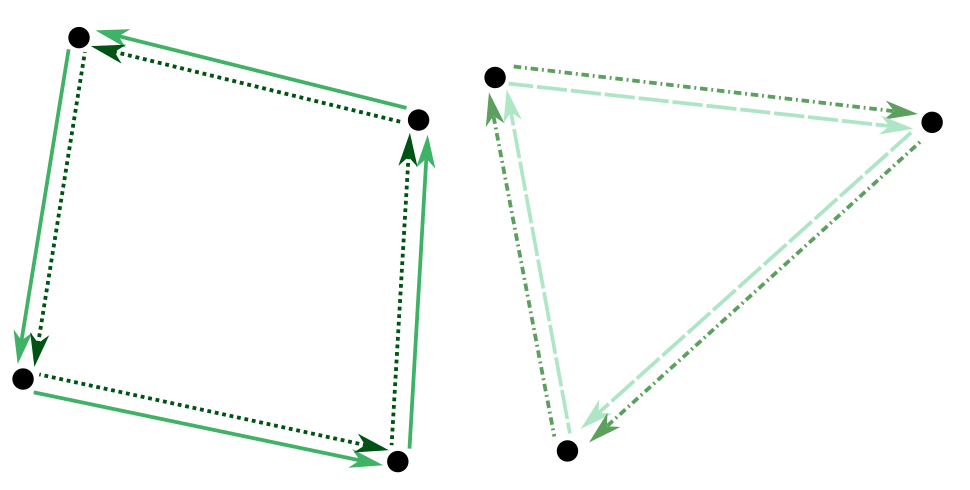}
  \caption{A graph formed by two $4$-directed cycles
    and two $3$-directed cycles satisfying
    the condition in 
    \eqref{eq:PartitionCondition}. 
    \label{fig:GraphOfCyclesPartition}}
\end{figure}

\noindent
Then, we suppose that $\mathrm{deg}(v) = 2$ for every $v \in V^\mu$. Choose
$(g_1,\dots,g_m)$ such that $[(g_1,\dots,g_m)]=\mu$. By using that all edges
of $G^\mu$ are multiple and that every vertex has degree exactly $2$ we can
see that there must be a partition into pairs of $\{1,\dots,m\}$ such that
(see Figure \ref{fig:GraphOfCyclesPartition})
\begin{equation}\label{eq:PartitionCondition}
  g_i = g_j\quad\text{if }i \sim j
  \quad\text{and}\quad  
  V^{g_i} \cap V^{g_j} = \emptyset \quad\text{if }i \not\sim j
\end{equation}
where the relation $\sim$ denotes if the elements belong to the same set of
the partition and $V^{g_l}$ denotes the vertex set of $g_l$ as before.
Necessarily $m$ is even and for each such $\mu$ we have
\[
  \frac{
    \mathrm{card}_n\left(\mu 
%      \cap 
%      \mathcal C_{k_1}^{(n)}
%      \times \dots \times 
%      \mathcal C_{k_m}^{(n)} 
    \right)}
  {n^{(k_1+\dots+k_m)/2}}
% \frac{1}{n^{(k_1+\dots + k_m)/2}}\frac{n^{(k_1+\dots + k_m)/2}}
%  {\sqrt k_1 \cdots \sqrt k_m}
  \xrightarrow[n \to \infty]{}
  \frac{1}{\sqrt k_1 \cdots \sqrt k_m}
\]
where the term $\sqrt k_1 \cdots \sqrt k_m$ appears because we are counting
cycles with no distinguished vertex. There is precisely one $\mu$ associated
to any partition $\mathcal P$ into pairs of $\{1,\dots,m\}$ such that
\begin{equation}\label{eq:PartitionCondition2}
  k_i = k_j \quad\text{if } i \sim j.
\end{equation}
Then, we shall use the notation $W_{\mathcal P} = W_\mu$ to notice that
\[
  \frac{1}{n^{(k_1+\cdots+k_m)/2}}\sum_{\mu\in\mathcal T_{(k_1,\ldots,k_m)}^{k_1+\cdots+k_m}}
  \mathrm{card}_n\left(\mu 
%    \cap 
%    \mathcal C_{k_1}^{(n)}
%    \times \cdots \times 
%    \mathcal C_{k_m}^{(n)} 
  \right)
  W_\mu
  \xrightarrow[n \to \infty]{}
  \frac{1}{\sqrt k_1 \cdots \sqrt k_n}
  \sum_{\mathcal P}
  W_{\mathcal P}
\]
where the sum is over all partition into pairs of $\{1,\dots,m\}$ such that
\eqref{eq:PartitionCondition2} happens. Since
\[
  W_{\mathcal P}= \prod\mathbb E \left[\alpha_{k_i}^{s_{i}}\alpha_{k_j}^{s_{j}} \right],
\]
where the product runs over all the pairs $\{i,j\}$ with $i \sim j$ and
$i \neq j$, and where $\alpha_{k} = a_{12}a_{23}\dots a_{(k-1)k}a_{k1}$. We
may use Isserlis/Wick theorem to conclude.

\subsection{Proof of Lemma \ref{le:meancv}}
\label{se:le:meancv}
\begin{figure}[htbp]\centering
  \includegraphics[width=0.3\columnwidth]{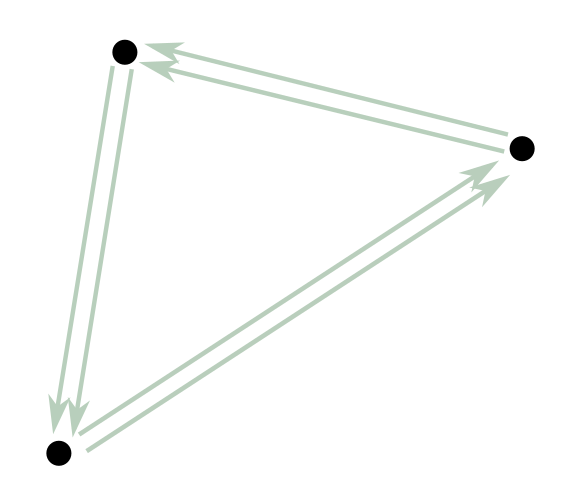}
  \caption{A single graph formed by a double cycle that counts for obtaining
    the expected value. \label{fig:DoubleCycle}}
\end{figure}

We start by checking that
\begin{equation}\label{eq:Erk}
  \mathbb E \left[ 
    \frac{r_k^{(n)}}{n^{k/2}}\right]
  \xrightarrow[n \to \infty]{} 
  \begin{cases}
    \mathbb E\left[{a_{11}}^2\right]^{k/2}
    & \text{if $k$ is even}\\
    0
    & \text{ if $k$ is odd.}
  \end{cases}
\end{equation}
We may use the same kind of counting argument as in Lemma \ref{le:clt}. Given
a sequence $(i_1,\dots,i_k)$, we construct a multigraph from it and notice
that every edge must be multiple for the graph to contribute. Next, if some
vertex has outer degree greater or equal than three then that graph does not
contribute neither. Finally, if the graph constructed from $(i_1,\dots,i_k)$
has every edge multiple and every vertex has outer degree two we can show that
it is a double cycle (see Figure \ref{fig:DoubleCycle}), in other words $k$ is
even, and for $l < l'$ we have that
\[
  i_l = i_{l'} \quad \text{ if and only if } \quad l' = l +
  \frac{k}{2}.
\]
The expectation $\mathbb E [a_g] $ for such double cycle $g$ is equal to $\mathbb E [a_{11} ^2]^{k/2}$. Since there are
\[
%	(k/2)!\binom{n}{k/2}
%´  =
  n(n-1)\dots\Bigr(n-\frac{k}{2}+1\Bigr)
\]
of those $(i_1,\dots,i_k)$, we have checked that \eqref{eq:Erk} holds. 

Now, for any pair of square-integrable
complex random variables 
$X$ and $Y$ let us use the notation
\[\mathrm{cov}(X,Y) = 
\mathbb E \big[ (X-\mathbb E[X])
(\overline Y-\mathbb E[\overline Y])
\big] \quad \mbox{ while }
\quad 
\mathrm{var}(X) = 
\mathrm{cov}(X,X)
=\mathbb{E}[|X-\mathbb{E}[X]|^2].
\]
To complete the proof of Lemma \ref{le:meancv}, it is sufficient to prove that
\begin{equation}\label{eq:varrk}
             \mathrm{var}\left( \frac{r_k^{(n)}}{n^{k/2}} \right) \xrightarrow[n \to \infty]{} 0.
             \end{equation}  
To this end, if $i = (i_1,\ldots,i_k)$, we set 
\[
a_i = a_{i_1 i_2}
  a_{i_2 i_3}\dots a_{i_{k-1} i_k}
  a_{i_k i_1}.
\]
%\[
%\widetilde a_i = a_{i_1 i_2}
%  a_{i_2 i_3}\dots a_{i_{k-1} i_k}
%  a_{i_k i_1} - \mathbb E [ a_{i_1 i_2}
%  a_{i_2 i_3}\dots a_{i_{k-1} i_k}
%  a_{i_k i_1}].
%\]
By construction, we have
%\begin{equation}\label{eq:varrkn1}
%\mathrm{var}\left(r_k^{(n)} \right)  = \sum_{i,j} \mathbb E \left[ \widetilde a_i ,
%\overline{\widetilde a_j}\right],
%\end{equation}
\begin{equation}\label{eq:varrkn1}
\mathrm{var}\left(r_k^{(n)} \right)  = \sum_{i,j} \mathrm{cov} \left(  a_i ,
 a_j\right),
\end{equation}
where the sum is over all pairs $(i,j)$ of $k$-tuples such that both $i$ and
$j$ have less than $k$ distinct elements.

From Cauchy\,--\,Schwarz inequality, the following crude bound holds:
\[\left|\mathrm{cov} \left(  a_i ,
 a_j\right)\right|\leq 4 M^{2k}, \]
%\[
%\mathbb E [ |\widetilde a_i 
%\overline{\widetilde a_j} | ] \leq 4 M^{2k},
%\]
where $M$ is such that the support of $a_{11}$ is contained in the ball of
radius $M$. Also, as above, we may identify each $k$-tuple $i$ with a path of
length $k$. Setting $i_{k+1} = i_1$, let us introduce the set of visited
vertices and the set of directed edges by
\[
  V_i = \{ i_1,\ldots,i_k\}
  \quad\text{and}\quad
  E_i = \{ (i_l,i_{l+1}) : l = 1,\ldots,k \}.
\]
Then $G_i = (V_i,E_i)$ is the directed graph associated to $i$ (self-loop
edges allowed). We define its excess as
\[
  \chi_i = \mathrm{card }(E_i) - \mathrm{card }(V_i) + 1 \geq 0.
\]
It is the minimal number of edges to be removed such that the remaining
subgraph has no undirected cycle (with the convention that $(u,u)$ is a cycle
of length $1$ and for $u \ne v$, $\{ (u,v),(v,u)\}$ forms a cycle of length
$2$). Since $G_i$ is the graph associated to a path of length $k$, the
assumption $\mathrm{card }(V_i) < k$ implies that
\[
  \chi_i \geq 2.
\]
Similarly, if $i,j$ are two $k$-tuples, we consider their associated graph
with vertex and directed edge sets
\[
  V_{ij} = V_i \cup V_j
  \quad\text{and}\quad
  E_{ij} = E_i \cup E_j.
\]
The excess of the corresponding graph $G_{ij} = (V_{ij},E_{ij})$ is
\[
  \chi_{ij} = \mathrm{card }(E_{ij}) - \mathrm{card }(V_{ij}) + c
\]
where $c \in \{ 1,2\}$ is the number of weak connected components of $G_{ij}$:
$c=1$ if $V_i \cap V_j \ne \emptyset$ and $c=2$ otherwise. Since $G_{ij}$ is
the union of $G_i$ and $G_j$, we have
\[
\chi_{ij} \geq \max ( \chi_i,\chi_j) \geq 2.
\]
Now, from the independence of the entries of the matrix $A_n$, we have
%$\mathbb E \left[ \widetilde a_i 
%\overline{\widetilde a_j}\right] = 0$ 
$\mathrm{cov}(a_i,a_j)=0$ unless $E_i\cap E_j$ is not empty. Thus $G_{ij}$ is
connected for such $i,j$. Moreover,
%$\mathbb E \left[ \widetilde a_i 
%\overline{\widetilde a_j}\right] = 0$ 
$\mathrm{cov}(a_i,a_j)=0$ unless all edges of $E_{ij}$ are visited at least
twice by the union of paths $i$ and $j$. Hence, for such $i,j$,
$\mathrm{card}(E_{ij}) \leq k$ and thus
\[
\mathrm{card}(V_{ij}) = 1 - \chi_{ij}  + \mathrm{card}(E_{ij})  \leq k-1.
\]
We thus have checked that 
\[
\mathrm{var}\left(r_k^{(n)} \right)  \leq 4C_k M^{2k} n^{k-1},
\]
where $C_k$ bounds the number of possibilities for the pair of $k$-tuples
$(i,j)$ once the set $V_{ij}$ is chosen. This gives \eqref{eq:varrk}, which
concludes the proof of the lemma.

\bibliographystyle{amsplain}
\bibliography{radiusbis}

\end{document}